\documentclass[12pt]{article}
\newtheorem{prop}{Proposition}
\newtheorem{corollary}{Corollary}
 \newtheorem{thm}{Theorem} 
\newtheorem{defn}{Definition} \input{diagrams}
\usepackage{amsmath} 
\usepackage{array} \newcolumntype{C}{>{$}c<{$}}
\begin{document}
\title{Pregroupoids and their enveloping groupoids} \author{Anders Kock}
\date{} \maketitle

\small ABSTRACT.  We prove that the forgetful functor from groupoids 
to pregroup\-oids has a left adjoint, with the front adjunction 
injective.  Thus we get an {\em enveloping groupoid} for any 
pregroupoid.  We prove that the category of torsors is equivalent to 
that of pregroupoids.  Hence we also get enveloping groupoids for 
torsors, and for principal fibre bundles.

\normalsize
\medskip

\section*{Introduction} The present note advocates the algebraic 
notion of {\em pregroupoid} as a natural context in which to study and compare  
groupoids, principal fibre bundles, torsors, bitorsors.  The aim has 
been to provide a theory which is {\em functorial}, and can 
immediately be interpreted in a wide variety of categories, in 
particular, in all toposes.  Hence, the construction principle ``choose a 
base point'' is not used, since it violates not only functorality, but 
also violates the choice principles available in toposes, where 
``non-empty'' (= inhabited) objects may have no ``points'' (= global 
sections).  One motivation I had for looking for such a theory, was to 
have an adequate, purely algebraic, framework for studying connections 
in principal fibre bundles, cf.\ \cite{APFBC}, in the context of 
synthetic differential geometry, where topos methods are crucial.

The main construction in this  framework is the 
construction of an ``enveloping groupoid'' $ X^+$ of a pregroupoid 
$X$.  It in fact provides a left adjoint for the forgetful functor 
from groupoids to pregroupoids, and the unit for the adjunction is 
injective, whence the choice of the adjective ``enveloping''.  In 
particular, the functor $X\mapsto X^+$ is faithful.  -- As a special 
case, the construction provides an enveloping groupoid of any principal 
fibre bundle.

The enveloping groupoid construction can be described without the
notion of pregroupoid; this was in fact done in \cite{APFBC} (for the
case of principal fibre bundles).

An essential ingredient in the construction of $X^+$ is the Ehresmann ``edge 
groupoid'' $XX^{-1}$ of a principal fibre bundle $X$.  This edge 
groupoid  construction was carried out in the context of 
pregroupoids in \cite{LLN} (but in a less equational manner).  The functor 
$X\mapsto XX^{-1}$, however, unlike $X\mapsto X^+$, is not an adjoint, 
and is not faithful.

\medskip

I want to acknowledge a heated but fruitful e-mail exchange with Ieke 
Moerdijk in the Summer and Fall of 2002,  on some of the 
topics of the present paper and  \cite{APFBC}.

\section{Equational theory of pregroupoids}

We consider a groupoid ${\bf G}= G_{1}\begin{picture}(22,12)
\put(3,2){\vector(1,0){15}}
\put(3,7){\vector(1,0){15}}
           \end{picture} G_{0}$.
For any two subsets $A\subseteq G_{0}$ and $B\subseteq G_{0}$, we
let ${\bf G}(A,B)$ denote the set of arrows $\in G_1$ whose domain is in $A$ and
whose codomain is in $B$. If $A=B$, this carries structure of
groupoid, the full subgroupoid on $A$, which thus here is denoted
${\bf G}(A,A)$.

There are evident book-keeping maps $$d_{0}:{\bf G}(A,B)\to A \; \mbox{ and 
}\; d_{1}:{\bf G}(A,B)\to B.$$

In ${\bf G}$, we compose from left to right, denoting composition by
$\circ$.  Then composition in ${\bf G}$, together with the book-keeping
maps, provide $X={\bf G}(A,B)$ with a certain partially defined algebraic
structure: a ternary operation denoted $yx^{-1}z$, defined whenever
$d_1 (x)= d_1 (y)$ and $d_0 (x)= d_0 (z)$, (and then $d_0 (yx^{-1}z) =
d_0 (y)$ and $d_1 (yx^{-1}z) = d_1 (z)$).  Namely $$yx^{-1}z := y\circ
x^{-1} \circ z.$$
The reader may find the following display useful.  The vertices are
elements of $A$ and $B$, respectively (with $A$-objects in the left
hand column, $B$-objects in the right hand column).

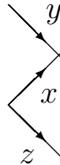
\begin{center} \begin{picture}(60,80) %%\put(60,80){\line(1,1){15}}
\put(10,80){\line(1,-1){20}}
\put(10,80){\vector(1,-1){14}}
\put(24,75){$y$}
\put(22,43){$x$}
\put(14,20){$z$}
\put(10,42){\line(1,1){20}}
\put(10,42){\vector(1,1){14}}
\put(10,42){\line(1,-1){20}}
\put(10,42){\vector(1,-1){14}}
\end{picture}
\end{center}

The following equations trivially hold for this ternary operation (whenever
the expressions are defined):
\begin{equation}xx^{-1}z = z \label{U1}\end{equation}
\begin{equation}yx^{-1}x = y \label{U2}\end{equation}
(``unit laws''),
      \begin{equation}vy^{-1}(yx^{-1}z) = vx^{-1}z \label{G1}\end{equation}
\begin{equation}(yx^{-1}z)z^{-1}w = yx^{-1}w, \label{G2}\end{equation}
(``concatenation laws''). The reason for the latter name is motivated
by the following diagrammatic device (also used in \cite{LLN}):

      We indicate  the assertion that $u=yx^{-1}z$ by a diagram

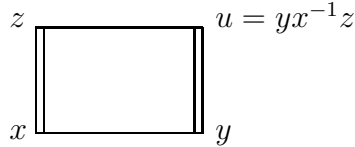
\begin{center}

\begin{picture}(60,70)
\put(10,10){\line(1,0){63}}
\put(10,50){\line(1,0){63}}
\put(10,10){\line(0,1){40}}
\put(13,10){\line(0,1){40}}
\put(70,10){\line(0,1){40}}
\put(73,10){\line(0,1){40}}
\put(0,8){$x$}
\put(78,8){$y$}
\put(0,50){$z$}
\put(78,50){$u = yx^{-1}z$}
\end{picture}

\end{center}

Here, the single lines connect elements in $X$ which
have same codomain, double lines connect elements with same domain.
Quadrangles that arise in this way, $u=yx^{-1}z$, we shall call {\em good} 
quadrangles, and (\ref{G1}) (resp.\ (\ref{G2})) then expresses that 
good quadrangles may be concatenated horizontally (resp.\ vertically).  
The display of (\ref{G1}) in terms of quadrangles in fact is

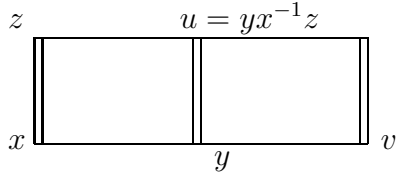
\begin{center}
\begin{picture}(120,70)
\put(10,10){\line(1,0){63}}
\put(10,50){\line(1,0){63}}
\put(10,10){\line(0,1){40}}
\put(13,10){\line(0,1){40}}
\put(70,10){\line(0,1){40}}
\put(73,10){\line(0,1){40}}
\put(136,10){\line(0,1){40}}
\put(133,10){\line(0,1){40}}
\put(141,8){$v$}
\put(0,8){$x$}
\put(78,2){$y$}
\put(0,54){$z$}
\put(65,54){$u = yx^{-1}z$}
\put(73,10){\line(1,0){63}}
\put(73,50){\line(1,0){63}}
\end{picture}
\end{center}

We proceed to make some purely equational deductions from  
(\ref{U1}), (\ref{U2}), (\ref{G1}), (\ref{G2}), which we take as 
axioms for the notion of pregroupoid.  To be specific, we pose

\begin{defn} A {\em pregroupoid} (``on $A,B$'') is an inhabited set $X$
equipped with surjections $\alpha :X\to A$, $\beta :X\to B$ and with a 
partially defined ternary operation, denoted $yx^{-1}z$, defined 
whenever $\beta (x)=\beta (y)$ and $\alpha (x)= \alpha (z)$; and then 
$\alpha (yx^{-1}z)= \alpha (y)$ and $\beta (yx^{-1}z) = \beta (z)$; and 
the equations (\ref{U1}), (\ref{U2}), (\ref{G1}), (\ref{G2}) are 
supposed to hold.
\end{defn}

(In \cite{LLN}, essentially the same notion was considered, but from a less
equational viewpoint).

Since the ``primitive'' operation $yx^{-1}z$ has three variables,
equations quick\-ly are equations in five or more variables, and
therefore it is convenient to denote the variables $x_1 , x_2 ,
\ldots$.  So the basic operation is $x_1 x_2 ^{-1} x_3$ (note that 
$x_1 $ corresponds to $y$, $x_2$ to $x$).  In fact, to make the 
notation even more lightweight, we may drop the ``$x$'' in $x_1$, 
$x_2$, ..  and instead just use the symbols $1,2,\ldots$.  So for 
instance, (\ref{G1}) is
\begin{equation} 41^{-1}(12^{-1}3) =42^{-1}3. \label{G1+}\end{equation}
-- The first
equational consequence of the axioms is an ``associative law'':

\begin{equation}(12^{-1}3)4^{-1}5 = 12^{-1}(34^{-1}5)
\label{ASS}\end{equation}
(provided the book-keeping makes the expressions meaningful, i.e.\
provided $\beta (1) = \beta (2), \alpha (2)=\alpha (3), \beta (3)= \beta 
(4), \alpha (4)=\alpha (5)$).

For,
$ (12^{-1}3)3^{-1}(34^{-1}5)$ equals $(12^{-1}3)4^{-1}5$, by (\ref{G1}),
and equals $12^{-1}(34^{-1}5)$ by (\ref{G2}).

Next, we have
\begin{equation} 21^{-1}(12^{-1}3)=3.\label{vertical}\end{equation}
For,  if we put
$4=2$ in (\ref{G1+}), we get $21^{-1}(12^{-1}3) = 22^{-1}3$ which is $3$,
by (\ref{U1}). Similarly, from (\ref{G2}) and (\ref{U2}), we get
\begin{equation}(12^{-1}3)3^{-1}2 =1.\label{horizontal}\end{equation}

        In the
diagrammatic form of quadrangles, as above, these two equations
express the following symmetry property for good quadrangles, which we
shall use without further comment in the ``graphical'' calculations
that follow.

\begin{prop}The mirror image og a good quadrangle in a horizontal
line, or in a vertical line, is again a good quadrangle. (So the
``Four-Group'' acts on the set of good quadrangles.)
\label{four-group} \end{prop}

{\bf Proof.}   Assume we have a good quadrangle
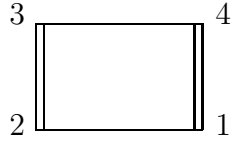
\begin{center}

\begin{picture}(60,70)
\put(10,10){\line(1,0){63}}
\put(10,50){\line(1,0){63}}
\put(10,10){\line(0,1){40}}
\put(13,10){\line(0,1){40}}
\put(70,10){\line(0,1){40}}
\put(73,10){\line(0,1){40}}
\put(0,8){$2$}
\put(78,8){$1$}
\put(0,50){$3$}
\put(78,50){$4$}
\end{picture}

\end{center}
so $4=12^{-1}3$. The fact that the reflection in a horizontal line is
a good quadrangle is the assertion that $1=43^{-1}2$, or, by the
assumption on $4$, that $1=(12^{-1}3)3^{-1}2$, which is just
(\ref{horizontal}). The assertion about reflection in a vertical line
similarly follows from (\ref{vertical}).

\medskip

We shall use the graphical calculus with good quadrangles to establish the 
following equation \begin{equation}6(34^{-1}5)^{-1}2 = 
65^{-1}(43^{-1}2);
\label{three-quad}\end{equation}
again the book-keeping conditions are assumed to make the expression
meaningful; these conditions are stated in diagrammatic form in the
diagram
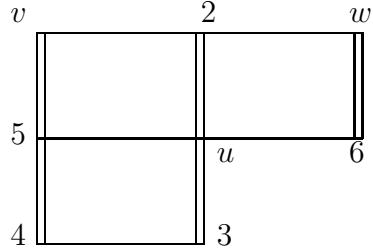
\begin{center}

\begin{picture}(115,90)(0,-35)
\put(10,10){\line(1,0){63}}
\put(10,50){\line(1,0){63}}
\put(10,10){\line(0,1){40}}
\put(13,10){\line(0,1){40}}
\put(70,10){\line(0,1){40}}
\put(73,10){\line(0,1){40}}
\put(70,10){\line(1,0){63}}
\put(70,50){\line(1,0){63}}
\put(130,10){\line(0,1){40}}
\put(133,10){\line(0,1){40}}
\put(10,-30){\line(1,0){63}}
\put(10,-30){\line(0,1){40}}
\put(13,-30){\line(0,1){40}}
\put(70,-30){\line(0,1){40}}
\put(73,-30){\line(0,1){40}}
\put(70,10){\line(0,1){40}}
\put(73,10){\line(0,1){40}}
\put(70,10){\line(0,1){40}}
\put(73,10){\line(0,1){40}}

\put(0,8){$5$}
\put(78,1){$u$}
\put(0,54){$v$}
    \put(72,54){$2$}
    \put(128,1){$6$}
    \put(128,54){$w$}
    \put(0,-30){$4$}
    \put(78,-30){$3$}
\end{picture}
\end{center}
Assume that all the three displayed quadrangles are good.  (These 
three quadrangles are constructed out of the data of the entries $2,3,4,5,6$, 
by first constructing $u$, and then $v$ and $w$.) So $u= 34^{-1}5$, and 
also $v=43^{-1}2$ (concatenate the two left hand quadrangles).  So by 
concatenating the two top quadrangles, we get
$$w = 65^{-1}v = 65^{-1}(43^{-1}2);$$
on the other hand, the upper right quadrangle expresses that $$w=6u^{-1}2 =
6(34^{-1}5)^{-1}2;$$
comparing, we get (\ref{three-quad}).

\medskip
Let us call a pair $(x,z)$ with $\alpha (x)= \alpha (z)$ a {\em
vertical} pair. The (horizontal) concatenation property for good
quadrangles, together with (\ref{U1}) and one of the symmetries 
mentioned in the Proposition, imply that we get an {\em equivalence} 
relation $\cong _v$ on the set of vertical pairs, namely
$$(x,z)\cong _v (y,u) \mbox{ iff } u=yx^{-1}z,$$
(geometrically: $(x,z)$ and $(y,u)$ form the vertical sides of a good
quadrangle). Note that if $(x,z)\cong _v (y,u)$, then $\beta (x)=
\beta (y)$ and $\beta (z)= \beta (u)$. The equivalence class of
$(x,z)$ is denoted $x^{-1}z$. The set of such equivalence classes is
denoted $X^{-1}X$.

Similarly, let us call a pair $(x,y)$ with $\beta (x)= \beta (y)$ a {\em
horizontal} pair. The (vertical) concatenation property for good
quadrangles, together with (\ref{U2}) and one of the symmetries 
mentioned in the Proposition, imply that we get an {\em equivalence} 
relation $\cong _h$ on the set of horizontal pairs, namely
$$(x,y)\cong _h (z,u) \mbox{ iff } u=yx^{-1}z,$$
(geometrically: $(x,y)$ and $(z,u)$ form the horizontal sides of a good
quadrangle). Note that if $(x,y)\cong _h (z,u)$, then $\alpha (x)=
\alpha (z)$ and $\alpha (y)= \alpha (u)$. The equivalence class of
$(x,y)$ is denoted $yx^{-1}$. The set of such equivalence classes is
denoted $XX^{-1}$.

\medskip
We proceed to derive some equations that involve ``fractions''
$yx^{-1} \in XX^{-1}$ and $x^{-1}z \in X^{-1}X$.  Among these is
\begin{equation}(12^{-1}3)4^{-1} =1(43^{-1}2)^{-1}.
\label{10}\end{equation}
By definition of the equivalence relation that defines $XX^{-1}$, this means
$$1=(12^{-1}3)4^{-1}(43^{-1}2).$$
But for the right hand side here, we have that it equals
$(12^{-1}3)3^{-1}2$ (using (\ref{G1})), which in turn by (\ref{G2}) equals
$12^{-1}2$, which is $1$, by (\ref{U2}).  -- Similarly, one proves
\begin{equation}(43^{-1}2)^{-1}5=2^{-1}(34^{-1}5)
\label{11}\end{equation}

\section{Enveloping groupoid of a pregroupoid}

Since the notion of pregroupoid is purely algebraic (except for the
surjectivity requirement for $\alpha $ and $\beta$),
it is clear how to organize pregroupoids into a
category (it will be upgraded into a 2-category in Section 4): if 
$$\begin{diagram} A& \lTo ^{\alpha}&X& \rTo ^{\beta}& B
      \end{diagram}$$
      and
       $$\begin{diagram} A'& \lTo ^{\alpha '}&X'& \rTo ^{\beta '}& B'
      \end{diagram}$$
      are pregroupoids, a morphism $\xi$ from the first to the second
      consists of maps $\xi _{0}:A\to A'$, $\xi _{1}: B\to B'$ and $\xi
       :X \to X'$ commuting with the structural maps and preserving
      the ternary operation. It is usually harmless to omit the subscripts
      and just write $\xi$ for all three maps in question.

In the following, ``groupoid'' means ``inhabited groupoid''. There is an 
evident functor
$$\mbox{groupoids } \to \mbox{ pregroupoids}$$
taking the groupoid ${\bf G}= G_{1}\begin{picture}(22,12)
\put(3,2){\vector(1,0){15}}
\put(3,7){\vector(1,0){15}}
           \end{picture} G_{0}$ to the pregroupoid ${\bf G}(G_0 ,G_0 )$, (so the
ternary operation $yx^{-1}z$ is given by $y\circ x^{-1}\circ z$).

\begin{thm}This functor has a left adjoint; the front adjunction is
injective.
        \end{thm}

     {\bf Proof/construction.} Given a pregroupoid $X = (A\leftarrow X \to
     B)$,
     we construct a groupoid $X^{+}$ whose set of objects is the disjoint
     sum of $A$ and $B$. Thus the set of arrows will  be a disjoint sum of
     four sets, $$X^{+}(A,A),  X^{+}(B,B),  X^{+}(A,B), X^{+}(B,A).$$ We
     first describe these sets with their book-keeping maps  $d_0$ and
$d_1$, and
     then we describe the composition $\circ$:

     $$\begin{array}{lll}X^{+}(A,A) := XX^{-1};& d_0 (yx^{-1}):= \alpha
     (y),& d_1 (yx^{-1}):=
     \alpha (x)\\
     X^{+}(B,B) :=X^{-1}X;& d_0 (x^{-1}z):=\beta (x),& d_1 (x^{-1}z):=
     \beta (z)\\
     X^{+}(A,B):=X;& d_0 (x):= \alpha (x),& d_1 (x): = \beta (x);
     \end{array}$$
     and finally $X^{+}(B,A)$ is to be another copy of $X$, which we will
     denote $X^{-1}$. An element $x$ of $X$ will be denoted $x^{-1}$ when
     considered in this copy of $X$. Thus we put
     $$\begin{array}{lll}X^{+}(B,A):= X^{-1};& d_0 (x^{-1}):= \beta
(x),& d_1 (x^{-1}):=
     \alpha (x)\end{array}.$$

     The fact that the book-keeping maps for $yx^{-1}$ and $x^{-1}z$ are
     well defined follows from the remarks at the end of Section 1.

Here is the description of the composition in the form of a
multiplication table.  We compose from left to right; the type of the 
left hand factor left hand factor is listed in the column on the left, 
the type of the right hand factor in the row on the top.

     \begin{center}\label{tal1}
\begin{tabular}{C|C|C|C|C|}

&A\to A&A\to B&B\to A&B\to B\\ \hline
A\to A&12^{-1}\circ 34^{-1}&12^{-1}\circ 3&&\\
&:=(12^{-1}3)4^{-1}&:=12^{-1}3&&\\ \hline
A\to B&&&1\circ 2^{-1} &1\circ 2^{-1}3\\
&&&:=12^{-1}&:=12^{-1}3\\ \hline
B\to A &3^{-1}\circ 21^{-1}&2^{-1}\circ 3&&\\
&:=(12^{-1}3)^{-1}&:=2^{-1}3&&\\ \hline
B\to B &&&3^{-1}2\circ 1^{-1}&2^{-1}3\circ 4^{-1}5\\
&&&:=(12^{-1}3)^{-1}&:=2^{-1}(34^{-1}5)\\ \hline
\end{tabular}
\end{center}
We proceed to check that the operation $\circ$ thus defined is
associative.  There are 16 cases to be considered, namely one for each
4-letter word in the letters $A$ and $B$.  All these cases follow
directly from the defining equations in the table, together with
equtions already derived for the ternary operation $yx^{-1}z$.  The
cases $AAAA$, $AAAB$, $AABB$, $ABBB$ and $BBBB$ were in effect dealt
with, from a different viewpoint,  in \cite{LLN}.  In terms of
equations, the
case $AAAB$, for instance, is proved as follows: $$(12^{-1}\circ
34^{-1})\circ 5 =
(12^{-1}3)4^{-1}\circ 5 =
(12^{-1}3)4^{-1}5,$$
using two definitions from the table, and $$12^{-1}\circ (34^{-1}\circ 5)
=12^{-1}\circ(34^{-1}5) = 12^{-1}(34^{-1}5),$$
likewise using two definitions from the table; but these two
expressions are equal by virtue of (\ref{ASS}).

The cases not dealt with in \cite{LLN} are those eleven 4-letter words that
involve the phrase $BA$.
We proceed with these (eleven) cases.

$$AABA:\;\; 12^{-1}\circ 3 \circ 4^{-1}=
\left\{ \begin{array}{ll}
(12^{-1}3)\circ 4^{-1}&= (12^{-1}3)4^{-1} \\
12^{-1}\circ 34^{-1}&= (12^{-1}3)4^{-1}
\end{array} \right.$$
where here, and in the following similar calculations,
the first line indicates the bracketing $(x\circ y)\circ z$, the
second the bracketing $x\circ (y\circ z)$.
$$ABAA:\; \; 1\circ 2^{-1} \circ 34^{-1}
= \left\{ \begin{array}{ll}
12^{-1}\circ 34^{-1} &= (12^{-1}3)4^{-1}\\
1\circ(43^{-1}2)^{-1} &= 1(43^{-1}2)^{-1} \end{array} \right.$$
and these are equal by (\ref{10}).
$$BAAA:\;\; 2^{-1}\circ 34^{-1} \circ 56^{-1} =
\left\{ \begin{array}{ll}
(43^{-1}2)^{-1}\circ 56^{-1}&=(65^{-1}(43^{-1}2))^{-1}\\
2^{-1} \circ (34^{-1}5)6^{-1}& =(6(34^{-1}5)^{-1}2)^{-1} \end{array} \right.$$
and these are equal by (\ref{three-quad}).

$$ABAB:\;\; 1\circ 2^{-1}\circ 3
= \left\{ \begin{array}{ll}
12^{-1}\circ 3 &= 12^{-1}3\\
     1\circ 2^{-1}3&= 12^{-1}3 \end{array} \right.$$
     (This is the crucial case for the comparison of the pregroupoid $X$
     and its enveloping groupoid $X^+$ !)

     $$ABBA:\; \; 1\circ 2^{-1}3\circ 4^{-1}
= \left\{ \begin{array}{ll}
(12^{-1}3)\circ 4^{-1} &= (12^{-1}3)4^{-1}\\
1\circ(43^{-1}2)^{-1} &= 1(43^{-1}2)^{-1} \end{array} \right.$$
and these are equal by (\ref{10}).

$$BAAB:\; \; 2^{-1}\circ 34^{-1}\circ 5
= \left\{ \begin{array}{ll}
(43^{-1}2)^{-1}\circ 5 & = (43^{-1}2)^{-1}5\\
2^{-1}\circ (34^{-1}5) &= 2^{-1}(34^{-1}5)\end{array} \right.$$
and these are equal by  (\ref{11}).

$$BABA:\; \; 2^{-1}\circ 3\circ 4^{-1}=\left\{ \begin{array}{ll}
2^{-1}3 \circ 4^{-1} & =(43^{-1}2)^{-1} \\
2^{-1}\circ 34^{-1}& =(43^{-1}2)^{-1}\end{array} \right.$$

$$BBAA:\;\;2^{-1}3\circ 4^{-1} \circ 56^{-1}=
\left\{ \begin{array}{ll}
    (43^{-1}2)^{-1}\circ 56^{-1}&= (65^{-1}(43^{-1}2))^{-1}\\
    2^{-1}3\circ (65^{-1}4)^{-1}&= ((65^{-1}4)3^{-1}2)^{-1} \end{array}
\right.$$
    and these are equal by (\ref{ASS}).

$$BABB:\;\; 2^{-1}\circ 3 \circ 4^{-1}5=
\left\{ \begin{array}{ll}
2^{-1}3\circ 4^{-1}5&= 2^{-1}(34^{-1}5)\\
2^{-1}\circ (34^{-1}5)&=2^{-1}(34^{-1}5)
\end{array} \right.$$

$$BBAB:\;\; 2^{-1}3\circ 4^{-1}\circ 5 =\left\{ \begin{array}{ll}
(43^{-1}2)^{-1}\circ 5&= (43^{-1}2)^{-1}5\\
2^{-1}3\circ 4^{-1}5&= 2^{-1}(34^{-1}5)
\end{array} \right.$$
and these are equal by (\ref{11}).

$$BBBA:\; \; 2^{-1}3 \circ 4^{-1}5 \circ 6^{-1}
= \left\{ \begin{array}{ll}
2^{-1}(34^{-1}5)\circ 6^{-1} & = (6(34^{-1}5)^{-1}2)^{-1}\\
2^{-1}3 \circ (65^{-1}4)^{-1}&= ((65^{-1}4)3^{-1}2)^{-1}\end{array} \right.$$
and these are equal by (\ref{three-quad}) and (\ref{ASS}).

These calculations prove that the composition $\circ$ is associative.
To prove the existence of units for $X^{+}$, we use that $\alpha
:X\to A$ and $\beta :X\to B$ were assumed surjective. For $a\in A
\subseteq A+B$, pick an element $x\in X$ with $\alpha (x)=a$.  Then 
$xx^{-1}$ will serve as a unit for the object $a$; for, we have
$$xx^{-1}\circ uv^{-1} = (xx^{-1}u)v^{-1} = uv^{-1},$$
using (\ref{U1}). Also,
$$yz^{-1}\circ xx^{-1} = (yz^{-1}x)x^{-1}=y(xx^{-1}z)^{-1},$$
by (\ref{10}) (put $3=4$). Now
    use (\ref{U1}) to get $yz^{-1}$ back. Finally $xx^{-1}\circ u =
xx^{-1}u = u$, by
(\ref{U1}) again.

Similarly, one proves that $x^{-1}x$ is a unit for
the object $b=\beta (x)\in B\subseteq A+B$.

Inverses are also almost tautologically present:
$yx^{-1}$ has $xy^{-1}$ as an inverse; for
$$yx^{-1}\circ xy^{-1} = (yx^{-1}x) y^{-1} =yy^{-1},$$
     using (\ref{U2}). Similarly, $z^{-1}x$ and $x^{-1}z$ are mutually
     inverse; and finally, $x$ as an arrow $a\to b$, has $x^{-1} :b\to a$ as
     inverse.

     This proves that $X^+$ is a groupoid.

     The set $X^+ (A,B)$ is by construction of $X^+$ just the given pregroupoid
     $X$, so that we have an injective mapping $\eta$ from $X$ to the
     (set of arrows of) $X^+$. Since $y\circ x^{-1}\circ z = yx^{-1}z$
     (cf.\ the case $ABAB$ in the above proof), $\eta $ is a morphism of
     pregroupoids into the underlying pregroupoid of $X^+$, so it just
     remains to check its universal property.  A
     pregroupoid homomorphism from $A\leftarrow X\to B$ into the underlying
     pregroupoid of a
     groupoid ${\bf G}=  G_{1}\begin{picture}(22,12)
\put(3,2){\vector(1,0){15}}
\put(3,7){\vector(1,0){15}}
           \end{picture} G_{0}$, consists of $\phi _0 : A\to G_0$,
           $\phi _1 :B\to G_0$, and a map $\phi  :X\to G_1$. The
           maps $\phi _0$ and $\phi _1$ together define a map
          $A+B \to G_0$, which is the object part of the desired
          functor $\overline{\phi}:X^+ \to {\bf G}$. The value of
          $\overline{\phi}$ on arrows is forced to be
          $\overline{\phi}(x)=x$ (since we require $\overline{\phi}$
          composed with $\eta$ to give $\phi$), and then the remaining
          three cases are also forced if we want $\overline{\phi}$ to be
          a functor:
	$$\overline{\phi}(x^{-1}):=(\phi (x))^{-1},
	\overline{\phi}(xz^{-1}):=\phi (x)\circ \phi (z)^{-1},
	\overline{\phi}(yx^{-1}):= \phi (y)\circ \phi (x)^{-1},$$
	and it is clear from the defining formulas (from the table) that
	$\overline{\phi}$ preserves composition, and also clearly identities.
	This proves the Theorem.
	
\medskip

Because $X$ {\em embeds} into the groupoid $X^+$, we propose the
name {\em enveloping groupoid of $X$} for it. It is analogous to the
enveloping associative algebra of a Lie algebra in the sense that all
equations concerning the ternary operation $yx^{-1}z$  can be checked
under the assumption that $yx^{-1}z$ is actually $y\circ x^{-1} \circ
z$ for the (associative) composition $\circ$ of a groupoid. We have
utilized
this principle in the calculations concerning connections in principal
fibre bundles, cf.\ \cite{APFBC}, (where the name {\em comprehensive
groupoid} was used for what we here call enveloping groupoid).

The set of objects of the enveloping groupoid $X^+$ contains the sets
$A$ and $B$ as
subsets. With respect to these subsets it has a certain property,
namely it is $A$-$B$-{\em transitive}; by this we mean that to every
$a\in A$ there
exists a $b\in B$ and an arrrow $a\to b$; and to every $b\in B$,
there exists an $a\in A$ and an arrow $a\to b$. The first assertion
is an immediate consequence of the surjectivity of $\alpha :X\to A$,
the second of the surjectivity of $\beta$.

\section{Pregroupoids and torsors}

If ${\bf G}$ is a groupoid with object set $A$, and $\alpha :X\to A$
is a map, there is a well known notion of (left) {\em action}
of ${\bf G}$ on $X$: if $g$ is an arrow in ${\bf G}$ and $x\in X$  satisfies
$\alpha (x) = d_{1}(g)$, then $g\cdot x \in X$ is defined and $\alpha
(g\cdot x )=d_{0} (g)$.  Unit and associative laws are assumed.
Then $X$ becomes the set of objects
of a groupoid, the {\em action} groupoid of the action; the arrows are pairs
$(g,x)$ with $d_1 (g) = \alpha (x)$.

There is an evident category of left groupoid actions: an object is a
pair consisting of a groupoid $ {\bf G} = (G_{1}\begin{picture}(22,12)
\put(3,2){\vector(1,0){15}}
\put(3,7){\vector(1,0){15}}
           \end{picture} A$) and a map $\alpha :X\to A$ on which ${\bf G}$ acts; 
           and a morphism $({\bf G}, X\to A)\to ({\bf G}',X'\to A')$ is a pair 
           consisting of a functor ${\bf G}\to {\bf G}'$ and a map $X\to X'$, 
           which is compatible with the structural maps $X\to A$ and $X'\to A'$, 
           and with the actions.

           The category of right groupoid actions  is defined similarly.

            Finally, there is a category of bi-actions: an object consists a
span
            of maps $$\begin{diagram} A& \lTo ^{\alpha}&X& \rTo ^{\beta}& B
      \end{diagram}$$ and a pair of groupoids ${\bf G}$ and ${\bf H}$ acting
      on the left and right on $X\to A$ and $X\to B$, respectively, and so
      that the two actions commute with each other ($A$ being the object set 
      of ${\bf G}$, $B$ the object set of ${\bf H}$).

The category of left {\em torsors} is a full subcategory of the category
of left groupoid actions.  We take the notion of torsor in the
generality which was given to it by Duskin \cite{Duskin}.
   We say that the an action by ${\bf G}$ on $X\to A$ makes $X$ into a {\em
${\bf
   G}$-torsor} if $X\to A$ is surjective, and the action groupoid is an
   equivalence relation.  Also, $X$, or equivalently $A$, is assumed
   to be inhabited (``non-empty'').  So for $x,y \in X$, there is {\em at most}
   one $g\in G_{1}$ with $g\cdot x =y$; such $g$, then, may be denoted
   $yx^{-1}$.

            Consider a left ${\bf G}$-torsor structure on $\alpha :X\to A$.
Since the action groupoid is
an equivalence relation on $X$, we may consider its quotient $\beta :
X\to B$. We say that $X\to A$ is a left ${\bf G}$-torsor {\em over}
$B$, or that $B$ is the {\em orbit set} of the left ${\bf G}$-torsor.  
We may write $B=X/{\bf G}$.

   The category of right torsors is defined similarly as a full
   subcategory of the category of right actions.  Finally, a bitorsor is
   a bi-action on a span $A\leftarrow X \to B$, which is a left torsor and
   a right torsor, and so that the structural map $X\to B$ is a quotient 
   map for the action groupoid (equivalence relation) for the left 
   action, and $X\to A$ is a quotient map for the action groupoid for the 
   right action.  The category of bitorsors is then defined as the full 
   subcategory of the category of biactions consisting of bitorsors.

   We denote the three categories thus described $$lTORS, \;  rTORS,\; \mbox{ and }
    \; lrTORS,$$ respectively.

If ${\bf K}$ is a groupoid, and $A$ and $B$ are inhabited subsets of
its object-set, then the set $X:={\bf K}(A,B)$ carries a left action by the
groupoid ${\bf K}(A,A)$, by precomposition, (with $\alpha = d_0 :X \to
A$ as book-keeping map). If ${\bf K}$ is $A$-$B$-transitive in the
sense described at the end of the previous section, this ${\bf
G}(A,A)$-action is in fact a torsor. For, $\alpha :X\to A$ is
surjective, by the $A$-$B$-transitivity assumption, and to any $x, y
\in {\bf K}(A,B)$
there exists at most one $g\in {\bf K}(A,A)$ with $g\cdot x = y$.
There exists such $g$ precisely when $x$ and $y$ have same codomain
(assumed to be in $B$), and so $B$ is the orbit set of this torsor
$X$ (surjectivity of $X\to B$ follows again by the $A$-$B$-transitivity).

Similarly, ${\bf K}(A,B)$ carries a right action, by post-composition,
by the groupoid ${\bf K}(B,B)$, and is in fact a right torsor, under
assumption of $A$-$B$-transitivity.  It is in fact a ${\bf K}(A,A)$-${\bf
K}(B,B)$-bitorsor.  For, the actions commute, by the associativity of 
composition in ${\bf K}$.

In \cite{LLN}, we sketched (Example p.\ 199) how a torsor gives rise
to a pregroupoid.  We shall extend this result by showing that this
construction is in fact a description of  an equivalence between the
category of pregroupoids (as described in Section 2), and each of the
three categories $lTORS$, $rTORS$ and $lrTORS$.

We first recall the passage (functor) from left torsors to pregroupoids.
Let ${\bf G}$ be a groupoid, acting on the left on $\alpha :X\to A$
($A$ being the set of objects of ${\bf G}$), and assume that it makes
$X$ into a torsor with orbit set $B$ (with quotient map denoted $\beta
:X\to B$).  The action is denoted by a dot.

If now $\beta (x) = \beta (y)$, they are in the same orbit for the
action, and since the action is free (by the torsor condition), there
is precisely one $g\in {\bf G}$ with $g\cdot x = y$. This $g$ may
therefore be denoted by a proper name, we call it $yx^{-1}$, so that
$yx^{-1}\in {\bf G}$ is characterized by
$$yx^{-1}\cdot x = y.$$
We note that $d_1 (yx^{-1})=\alpha (x)$, and $d_0 (yx^{-1}) = \alpha
(y)$. If now $z\in X$ has $\alpha
(z)=\alpha (x)$, $yx^{-1}\cdot z$ makes sense, and $\alpha (yx^{-1}\cdot
z) = d_0 (yx^{-1}) = \alpha (y)$. Also, by construction, $z$ and
$yx^{-1}\cdot z$ are in the same orbit for the action, so $\beta (z)=\beta
(yx^{-1}\cdot z)$.  Thus, if we put
$$yx^{-1}z := yx^{-1}\cdot z,$$
we have defined a ternary operation on $X$ with the correct
book-keeping for a pregroupoid. We proceed to check the four
equations (\ref{U1})-(\ref{G2}). First, $xx^{-1}$  is the
identity arrow at $\alpha (x)$, by the unitary law of the ${\bf G}$-action.
So
$xx^{-1}\cdot z = z$, by the unitary law of the
${\bf G}$-action, hence $xx^{-1}z=z$, proving (\ref{U1}). The
equation (\ref{U2}), on the other hand, is just the defining equation for 
$yx^{-1}$.

Next consider $x,y,z,v \in X$ with $\beta (x)=\beta (y) = \beta (v)$
and $\alpha (x)=\alpha (z)$ (the reader may want to refer to the
graphic display -- a double quadrangle -- of precisely
these book-keeping conditions, in Section 1 above).  Then
$$vy^{-1}(yx^{-1}z)=vy^{-1}\cdot (yx^{-1}\cdot z) = (vy^{-1}\circ
yx^{-1})\cdot z,$$
using the defining equations (twice) for the first equality sign, and
the assiociative law for the ${\bf G}$-action for the second
(composition in ${\bf G}$ denoted by $\circ$). On the other hand
$$vx^{-1}z = vx^{-1}\cdot z$$
by definition, so  (\ref{G1}) will follow by proving
$$vy^{-1}\circ yx^{-1} = vx^{-1}.$$
By the freeness of the action, it suffices to prove
$$(vy^{-1}\circ yx^{-1})\cdot x = vx^{-1}\cdot x.$$
The left hand side is by the associative law of the action the same as
$vy^{-1}\cdot (yx^{-1}\cdot x)$, which is $v$, by two applications of
the defining equations for the ``fractions'' of the form $yx^{-1}$;
the right hand side is also $v$, by one such defining equation. This
proves (\ref{G1}).

Finally, consider $x,y,z,w \in X$ with $\beta (x)=\beta (y)$ and
$\alpha (x)=\alpha (z) = \alpha (w)$.  First, we claim that
$(yx^{-1}z)z^{-1} = yx^{-1}$; since the action is free, it suffices to
see that $(yx^{-1}z)z^{-1}\cdot z = yx^{-1}\cdot z$.  But both sides
of this equation are equal to $yx^{-1}z$.  So
$$(yx^{-1}z)z^{-1}w = ((yx^{-1}z)z^{-1})\cdot w = yx^{-1}\cdot w = yx^{-1}w,$$
proving (\ref{G2}).  Note that no equational assumptions
(unitary or associative law of action) were used in the proof of
(\ref{U2}) and (\ref{G2}).

\medskip

It is clear that a morphism of left torsors, in the sense explained
above, gives rise to a morphism of pregroupoids, so that the
construction described is actually a functor
\begin{equation}c_l : lTORS \to 
\mbox{pregroupoids}.\label{eq1}\end{equation} (For a left torsor $X\to 
A$ with orbit set $X\to B$, the pregroupoid constructed is a 
pregroupoid on $A,B$.)

Similarly, we have a functor $c_r : rTORS \to \mbox{pregroupoids}$.

On the other hand, the envelope construction provides a functor 
$$\mbox{pregroupoids}\to \mbox{groupoids},$$ $X\mapsto X^+$.  If 
$A\leftarrow X\to B$ is a pregroupoid, then $X^+$ is a groupoid with 
object set $A+B$, and so we may form the $X^+ (A,A)$-$X^+ 
(B,B)$-bitorsor $X^+ (A,B)$ by the recipe in the beginning of the 
present Section.  So the construction of enveloping groupoid gives 
also gives rise to a functor 
\begin{diagram}\mbox{pregroupoids}& \rTo ^{env}& lrTORS .
\end{diagram}
Finally, there are the two obvious forgetful functors
\begin{diagram} lrTORS &\rTo^{U _r}& lTORS&&lrTORS &\rTo^{U _l}& rTORS
\end{diagram}
forgetting the right and left action, respectively.

We collect the functors described here together in the diagram
\begin{diagram}
lrTORS & \rTo ^{U_l}&rTORS&&\\
\dTo^{U_r} && \dTo_{c_r}&&\\
lTORS&\rTo _{c_l}&\mbox{pregroupoids}&\rTo _{env}&lrTORS.
\end{diagram}

\begin{thm}All functors exhibited here are equivalences.  Any 
endofunctor composed of functors (=cyclic composite) exhibited here is 
isomorphic to the relevant identity functor.  The square commutes on 
the nose.  The cyclic composite (\ref{cyclic}) (below) is, on the nose, 
the identity functor on the category of pregroupoids.
\label{equ}\end{thm}

In the next Section, we shall describe yet another category equivalent to 
these, namely the category of (groupoid-) fibrations over the groupoid 
${\bf I}$ (=the generic invertible arrow). -- Note that among the four categories 
proved equivalent in the 
Theorem, the category of pregroupoids is the most ``compact'', in the sense 
of involving least data; this is why it is possible to have certain 
{\em strict} equalities between functors with values in the 
category of pregroupoids.

\medskip

{\bf Proof.} We first prove that the square commutes.  Let $A\leftarrow X\to 
B$ be a bitorsor for left and right actions of the groupoids ${\bf G}= 
G\begin{picture}(22,12) \put(3,2){\vector(1,0){15}} 
\put(3,7){\vector(1,0){15}} \end{picture} A$ and ${\bf H}= 
H\begin{picture}(22,12) \put(3,2){\vector(1,0){15}} 
\put(3,7){\vector(1,0){15}} \end{picture} B$, respectively.  The two 
pregroupoids constructed by the two functors $c_l \circ U_r$ and $c_t 
\circ U_l$ both have $A\rightarrow X\to B$ for its underlying sets, so 
it suffices to see that the two ternary operations on $X$ agree.  
Consider $x,y,z$ satisfying the relevant book-keeping conditions for 
formation of the two possible $yx^{-1}z$.  So there are (unique) 
arrows $g\in G$ and $h\in H$ so that $g\cdot x =y$ and $x\cdot h = z$.  
Then $g\cdot z = (yx^{-1}z)_l$ and $y\cdot h = (yx^{-1}z)_r$ (with 
$(yx^{-1}z)_l$, resp.\ $(yx^{-1}z)_r$, denoting the ternary operation 
coming from the left, respectively right, torsor structure).  We then 
have $$(yx^{-1}z)_l = g\cdot z = g\cdot (x\cdot h) = (g\cdot x)\cdot h 
=
y\cdot h = (yx^{-1}z)_r ,$$
using for the middle equality sign that the two actions commute with each other.

We next prove that the composite
\begin{equation}\begin{diagram}\mbox{pregroupoids} &\rTo^{env} & 
lrTORS&\rTo^{U_r}&lTORS &\rTo ^{c_l}&\mbox{pregroupoids} 
\end{diagram}\label{cyclic}\end{equation} 
is the identity functor (on 
the nose).  Starting with a pregroupoid $A\leftarrow X \to B$, the 
composite of the two first functors here gives the left $X^+ 
(A,A)$-torsor $X^+ (A,B)=X$, with ternary operation given in terms of 
the composition $\circ$ in $X^+$.  But the composition of arrows from 
$A$ to $A$ with arrows from $A$ to $B$ in $X^+$ are precisely defined 
by the ternary operation in $X$, cf.\ the entry with address $(1,2)$ 
in the table which defines the composition $\circ$.

Next, we prove that the composite
\begin{equation}\begin{diagram} lTORS & \rTo ^{c_l} &\mbox{pregroupoids} &\rTo 
^{env} &lrTORS&\rTo^{U_r}&lTORS
\end{diagram}\end{equation}
 is isomorphic to the identity 
functor on $lTORS$ (a similar statement holds for $rTORS$).  Given a 
left ${\bf G}$-torsor $A\leftarrow X$, with $X\to B$ as orbit set, the 
underlying object of the ``new'' torsor is again $A\leftarrow X$, so we 
just have to provide isomorphisms between the acting groupoids, in 
this case ${\bf G} \cong X^+ (A,A)$ (compatible with the actions).  
The object sets of both ${\bf G}$ and $X^+ (A,A)$ are $A$.  The 
isomorphism on the arrow sets is given by sending $g:a\to a'$ into 
$yx^{-1}$ where $y = g\cdot x$, $x\in X$ any element of $X$ over the 
codomain of the arrow $g$.  The passage the other way takes a 
``fraction'' $yx^{-1}$ to the unique arrow $g$ with $g\cdot x =y$ (as 
also anticipated by the notation $yx^{-1}$ which we used for this $g$ 
in the discussion of torsors).

The same argument, applied twice (once on the left and once on the 
right) proves that the composite
\begin{diagram}
lrTORS &\rTo ^{U_r}&lTORS & \rTo ^{c_l} &\mbox{pregroupoids} &\rTo ^{env} 
&lrTORS
\end{diagram}
is isomorphic to the identity functor on $lrTORS$.
So the three functors displayed here provides a ``cycle'' of three 
arrows, with all three cyclic composites isomorphic to the identity 
functor of the respective vertex.  So all three of them are 
equivalences.  Similarly, also $c_r$ and $U_l$ are equivalences.  This 
proves the Theorem.

\medskip

\begin{corollary} There is an adjoint equivalence
\begin{diagram}lTORS&\rTo ^{ad} &rTORS
\end{diagram}
\end{corollary}

Namely, take $ad$ to be the composite
\begin{diagram}
lTORS & \rTo ^{c_l} &\mbox{pregroupoids} &\rTo ^{env} &lrTORS&\rTo^{U_l}&rTORS
\end{diagram}
The quasi-inverse is constructed similarly (replace all $l$'s by $r$'s 
and conversely).  It may also be denoted $ad$.  -- This 
$ad$-equivalence is classical, and plays (at least in the case of 
bitorsors over groups) an important role in Giraud's book, 
\cite{Giraud} III.1, where the notation $ad$ also appears.

 \medskip

A torsor (right, say) $({\bf G}, X\to B)$, where $B=1$ (and thus ${\bf 
G}$ is just a {\em group}) is usually called a {\em principal ${\bf 
G}$-bundle over $A$}, where $A$ is the orbit set.  The construction of 
the adjoint groupoid $ad (X)$, or gauge groupoid $XX^{-1}$, which is a 
groupoid with $A$ as object set, is classical, due to Ehresmann, 
\cite{Ehr}.  In our context, it appears as a full subgroupoid of 
$X^+$, which in this case is a groupoid with $A+B = A+1$ as object 
set.  Note that the functor $X\mapsto X^+$ is faithful (since $X^+$ 
contains $X$ as a subset), whereas $X \to XX^{-1}$ is not.  A 
description of the enveloping (=comprehensive) groupoid $X^+$ for the 
case of principal bundles was given in \cite{APFBC}.  The construction 
there was carried out without the notion of pregroupoid; but then the 
naturality and symmetry of the construction is not so visible.

\section{Fibrations over {\bf I}}

We discuss fibered categories ${\bf E}\to {\bf B}$,
see e.g.\ \cite{Giraud} for this notion.  It is well known that the 
fibres of such a fibration are groupoids if and only if all arrows in ${\bf 
E}$ are cartesian.  For fixed base category ${\bf B}$, we thus get the 
category of such ``fibrations-in-groupoids'' over ${\bf B}$.  If ${\bf B}$ 
happens to be itself a groupoid, then the total category ${\bf E}$ of a 
fibration-in-groupoids is also a groupoid.  We let ${\bf I}$ denote the 
groupoid containing the ``generic invertible arrow'', in other words, 
${\bf I}$ has two objects $a_0$ and $b_0$, and besides the two identity 
arrows, it has one arrow $i:a_0\to b_0$ and one arrow $i^{-1}:b_0\to 
a_0$, and no other arrows.  It can be described in very many ways; for 
instance, it is the enveloping groupoid of the terminal pregroupoid 
${\bf 1}=1\leftarrow 1\to 1$.

The following result is an application of the enveloping groupoid.

\begin{thm}The category of bitorsors (hence also the category of 
pregroup\-oids, by Theorem \ref{equ}) is equivalent to the category of 
fibrations-in-groupoids over ${\bf I}$ with inhabited total category.
\label{cylinder}\end{thm}

{\bf Proof.} Given a fibration $\gamma :{\bf X}\to {\bf I}$. Let $A$ be the
set of objects $a$ in ${\bf X}$ with $\gamma (a)=a_0$, and $B$ the set
of objects in ${\bf X}$ with $\gamma (b)=b_0$. If $X$ is inhabited,
then so are both $A$ and $B$.  Then ${\bf X}$ is $A$-$B$ transitive, in the 
sense of Section 2 (end).  For, given $b\in B$, take a (cartesian) arrow 
over $i$ with codomain $b$; it will be an arrow from an object in $A$ 
to $b$.  Similarly for a given object $a\in A$ (utilize $i^{-1}$).  
Therefore, ${\bf X}(A,B)$ is a ${\bf X}(A,A)$-${\bf X}(B,B)$-bitorsor.  
Conversely, given a ${\bf G}$-${\bf H} $ bitorsor $X$ (where the 
object sets of ${\bf G}$ and ${\bf H}$ are $A$ and $B$, respectively).  
Consider it as a pregroupoid $X$ on $A$, $B$ via the functor from 
bitorsors to pregroup\-oids, described in Theorem \ref{equ}.  Its 
enveloping groupoid $X^+$ is a groupoid with object set $A+B$.  We get 
a functor $\gamma : X^+ \to {\bf I}$, easily described ad hoc (mapping 
each $a\in A$ to $a_0$ etc.); alternatively apply the (left adjoint) 
functor $(-)^+ :$ pregroupoids $\to$ groupoids to the unique pregroupoid 
morphism $X\to {\bf 1}$.

\medskip
We note that if ${\bf X}\to {\bf I}$ is an inhabited fibration, then 
the inclusion of either of the two fibres, i.e.\ the ``end'' groupoids ${\bf 
X}(A,A)$ and ${\bf X}(B,B)$, is an equivalence of categories.  For, 
they are clearly full and faithful, and essential surjectivity follows 
from the $A$-$B$ transitivity. 

However, the {\em functor} which to a fibration ${\bf X}\to {\bf I}$ 
associates the groupoid in either end, say  ${\bf 
X}(A,A)$ is not an equivalence; it is not even faithful.  

\medskip

Now the category of inhabited fibrations-in-groupoids over ${\bf I}$
is in an evident way a 2-category; the 2-cells are just natural
transformations. Thus, the components of the 2-cells
(natural transformations) are
 vertical arrows.  In particular, 2-cells are invertible.  Since the 
 inclusions of each of the two end-groupoids (or edge groupoids) are 
 equivalences, it follows that a 2-cell between two functors over ${\bf 
 I}$ is completely given by its components on the objects in 
 the $a_0$-end, or by its components on the objects in the $b_0$-end.  In 
 particular, enriching the category of principal bundles (over varying 
 groups) into a 2-category only amounts to considering the category of 
 groups as a 2-category in the standard way (2-cells being given by 
 ``conjugation by an element in the codomain group'').

 From the 
 equivalence of the Theorem  follows that there is a 2-dimensional 
 structure (with all 2-cells invertible) on the category of 
 pregroupoids, and the rest of this section just consists in making 
 this 2-dimensi\-onal structure explicit.

So consider two inhabited fibrations ${\bf X}\to {\bf I}$ and ${\bf
X'}\to {\bf I}$, and two functors $f$ and $g: {\bf X}\to {\bf X'}$ over
${\bf I}$, and let $\tau :f \to g$ be a natural transformation.
Denote by $A$ and $B$ the set of objects in ${\bf X}$ over $a_0$ and 
$b_0$, respectively, and similarly $A'$ and $B'$ in ${\bf X'}$.  The 
functor $f$, being a functor over ${\bf I}$, induces maps $A \to A'$ 
and $B\to B'$, these maps are also just denoted $f$.  Similarly the 
maps $A \to A'$ and $B\to B'$ induced by $g$ are denoted $g$.  
Finally, let $X$ be the pregroupoid ${\bf X}(A,B)$ on $A,B$, and 
similarly for $X'$ on $A', B'$.

Consider for $a\in A$ the arrow $\tau _a :
f(a) \to g(a)$ in ${\bf X}'(A,A)$. For any $u:g(a)
\to b$ ($b$ an object $\in B'$; such $u$ exist by $A'$-
$B'$-transitivity of ${\bf X}'$), let $t(a,u): f(a) \to b$ denote
$\tau _a \circ u$. Then of course $\tau _a = t(a,u) \circ u^{-1}$.
Note that $u\in X'$.
Similarly, for $b \in B$, $\tau _b :f(b) \to g(b)$ may be written
$\tau _b =v^{-1}\circ s(b,v)$, where $s(b,v)$ is $v\circ \tau
_b$. Note that $v\in X'$.

The reader may find the following display helpful:
\begin{diagram}f(a)&&&&a&\rTo ^v & f(b)\\
\dTo ^{\tau _a}&\rdTo ^{t(a,u)}&&&&\rdTo _{s(b,v)}&\dTo _{\tau _b}\\
g(a)&\rTo _u &b&&&&g(b)
\end{diagram}

Now $t$ and $s$ are (partially defined) maps which satisfy three
equations, and together, encode the information of the 2-cell $\tau$
in pure pregroupoid terms. Precisely, $t(a,u)$ is defined whenever
$\alpha '(u)= g(a)$ ($\alpha '$ denoting domain formation $X' \to
A'$); and then $\alpha ' (t(a,u))=f(a)$, $\beta ' (t(a,u))=\beta '(u)$,
where $\beta ':X' \to B'$ is codomain formation. Similarly, $s(b,v)$
is defined whenever $\beta '(v) =f(b)$, and then $\beta '(s(b,v))=g(b), \alpha ' 
(s(b,v))=\alpha '(v)$.  The following equations hold (assuming that 
the book-keeping conditions make them meaningful); we omit the sign 
$\circ$ for composition in ${\bf X}$ and ${\bf X}'$: 
\begin{eqnarray}t(a,g(x))& =&s(b,f(x))\label{fifteen}\\
t(a,u)v^{-1}w&=&t(a,uv^{-1}w)\label{sixteen}\\
wu^{-1}s(b,v)&=&s(b,wu^{-1}v)\label{seventeen}
\end{eqnarray}
The equation (\ref{fifteen}) just follows from the naturality of $\tau$ with 
respect to $x:a\to b$.  For, consider the commutative naturality 
square (expressing naturality of $\tau$ with respect to $x:a\to b$) 
\begin{equation}\begin{diagram}f(a)&\rTo ^{f(x)}&f(b)\\
\dTo^{\tau _a}&\rdTo&\dTo_{\tau _b}\\
g(a)&\rTo _{g(x)}&g(b)
\end{diagram}\label{nat}\end{equation} 
The diagonal makes both triangles 
commute, and their commutativity express that the diagonal is, 
respectively, $t(a,g(x))$ and $s(b, f(x))$, which thus are equal.

For the equation (\ref{sixteen}), both sides are equal $\tau _a \circ u 
\circ v ^{-1} \circ w$, and for (\ref{seventeen}), both sides are equal to 
$w\circ u^{-1} \circ v\circ \tau _b$.

We now show that the data of such $t$ and $s$, satisfying the three
equations, come from a unique natural transformation $\tau$.

To define $\tau _a$ for $a\in A$, pick by $A$-$B$-transitivity an 
arrow $u:g(a)\to b$ and put
$$\tau _a := t(a,u)\circ u^{-1}.$$
That this is independent of the choice of $u$ is an immediate 
consequence of (\ref{sixteen}).  Similarly,
$$\tau _b := v^{-1}\circ s(b,v)$$
for some $v:a\to f(b)$; this is independent of choice of $v$ by 
(\ref{seventeen}).  It remains to check naturality of the $\tau$ thus 
constructed.  Now, $\tau$ is natural with respect to arrows $a\to b$ 
(for $a\in A, b\in B$); this follows from (\ref{fifteen}), by chosing 
$u:=g(x), v:= f(x)$ in the defining equations for $\tau _a$ and $\tau 
_b$, respectively (contemplate (\ref{nat}), now with the two 
expressions in (\ref{fifteen}) as diagonal).

But arrows of the form $a\to b$ (for $a\in A, b\in B$) generate ${\bf 
X}$ as a groupoid, so therefore, naturality of $\tau$ with respect to 
such arrows implies naturality with respect to all arrows in ${\bf X}$.

  \section{Examples} Let $A$ and $B$ be two smooth manifolds of 
  dimensions $n$ and $k$, say, with $n\geq k$, and consider a geometric 
  distribution $D$ on $A$ of codimension equal to the dimension $k$ of 
  $B$.  Let $X$ be the set of all 1-jets of maps from $A$ to $B$ with 
  $D$ as kernel.  Precisely, for each $a\in A$, consider the set $X_a$ 
  of 1-jets at $a$ of maps $f:A\to B$ such that the kernel of $df _a 
  :T_a (A) \to T_{f(a)}(B)$ is the linear subspace $D_a \subseteq T_a 
  (A)$.  For dimension reasons, then, $df _a$ is surjective.  Let $X$ be 
  the disjoint union of all the $X_a$'s.  Then $X$ is born with a map 
  $\alpha :X\to A$, but is also comes with a map $\beta :X\to B$, namely 
  to the 1-jet of $f$ at $a$, asssociate $f(a)\in B$.  (Actually $X$ is 
  a submanifold of the standard jet manifold $J^1 (A,B)$ of 1-jets of 
  maps from $A$ to $B$.)
  
  We shall equip this $A\leftarrow X \to B$ with a ternary operation 
  making it into a pregroupoid.  So let $x$, $y$ and $z$ be 1-jets with 
  $D$ as kernel, in the sense explained, represented by functions $f,g$ 
  and $h$.  Assume $\alpha (x) = \alpha (z), =a$, say, and $\beta 
  (x)=\beta (y), = b$, say.  Since $d_a f$ and $d_a h$ are surjective 
  linear maps with the same kernel $D_a$, there is a unique bijective 
  linear map $ \kappa: T_{f(a)}B \to T_{h(a)}B$ with $ d_a f\circ \kappa 
  = d_a h$ (composing from left to right).  By the Inverse Function 
  Theorem, there is locally around $f(a)$ a smooth map $k$ with 
  $d_{f(a)}k = \kappa$.  We put $yx^{-1}z$ equal to the 1-jet af $a' = 
  \alpha (y)$ of the composite $g\circ k$.  This makes sense, since 
  $g(a' ) = f(a)$ by the book-keeping assumption $\beta (y) = \beta 
  (x)$.
  
  The verification of the four equations is straightforward.  -- The 
  edge groupoids of this pregroupoid are the following: $X^{-1}X$ is the 
  groupoid of all invertible 1-jets $b\to b'$ from $B$ to itself; 
  $XX^{-1}$ is the groupoid of 1-jets $a\to a'$ from $A$ to itself which 
  ``take $D$ into $D$'', i.e.  1-jets at $a$ of functions $F$ such that $dF_a 
  :T_a (A) \to T_{a'} (A)$ maps $D_a$ into $D_{a'}$.
  
 This latter groupoid also occurs as edge groupoid of a principal $GL(k)$ 
 bundle $Y$ over $M$, namely the bundle of surjective linear maps $T_a A 
 \to {\bf R}^k$ with $D$ as kernel.  But note that there is no {\em 
 natural} way of mapping the pregroupoid $X$ to the pregroupoid $Y$; in 
 fact, such a map would amount to a framing of the tangent bundle of $B$.
 
 \medskip
 
 Pregroupoids $A\leftarrow X\to B$ with both $A$ and $B$ equal to the 
 1-point set were considered in \cite{ATMF} under the name {\em 
 pregroups}.  The two edge groupoids are in this case just groups; in 
 fact two groups which are isomorphic, but not canonically isomorphic, 
 unless they are abelian.  {\em Picking} an element in $X$ will provide 
 a specific isomorphism between the two edge groups.
 
 My contention is that the notion of pregroup is simpler  than that 
 of group.  In fact, in some cases, it precedes the notion of group in the 
 process of understanding.  How many ways can you put three pigeons 
 into three pigeon holes ?  Without knowing anything about neither 
 pigeons nor mathematics, most people will, after a moments reflection, 
 be able to answer ``six''.  This {\em number} of ways (or this {\em 
 set} of ways, as mathematicians prefer to say) carries canonically the 
 structure of pregroup (being a set of bijections from one set to 
 another), but does not carry structure of group.  The two edge 
 groupoids are of course both ``the'' symmetric group $S_3$, namely the 
 group of permutations of the given three pigeons, respectively of the 
 three given pigeon holes.  What is a permutation of three pigeons ?  
 ``Put the white pigeon in the place where the grey pigeon was, and put 
 the grey pigeon in the place where \ldots ".  Not a very natural thing 
 to do, and in any case is equivalent to describing permutations in 
 terms of the {\em places} the pigeons were occupying, before and after 
 the permutation.  These places may as well be called ``pigeon holes'', 
 and then we are precisely describing the elements of the 
 pigeon-permutation group in terms of fractions $yx^{-1}$ made out of 
 the pregroup.
 
 A more mathematical version of this comment is the following: What is ``the'' 
 symmetric group in three letters ?  What is the sense of the definite 
 article ``the''?  The group of permutations of the three letters $A, 
 a$, and $\alpha$ is not the {\em same} as the group of permutations of the 
 three letters $b,c, d$; these groups are not even {\em canonically isomorphic} 
 (which in mathematics is sufficient justification for using the 
 definite article).  For, an isomorphism between them depends on 
 {\em choosing} a bijection between the two three-letter sets; a 
 different choice may change the constructed isomorphism by a 
 conjugation.  This means that ``the symmetric group in three letters'' 
 is well defined only in the category of groups and conjugation classes 
 of group homomorphisms, i.e., ``{\em the} symmetric group in three 
 letters'' is an object in the category of ``liens'', or ``bands'', in the 
 termonology of \cite{Giraud} resp.\ \cite{Duskin}.

\end{document}